\documentclass[a4paper]{article} 

\usepackage{amsfonts,amsmath,amssymb,amscd}
\usepackage{latexsym}
\usepackage[french]{babel}
\usepackage[T1]{fontenc}

\newcommand{\coker}{\text{coker}\ }
\newcommand{\ind}{\text{Ind }}
\newcommand{\di}{D_A}
\newcommand{\un}{\text{U(1)}}
\newcommand{\U}{\text{U(2)}}

\newcommand{\SO}{\text{SO(4)}}
\newcommand{\spi}{\text{Spin}}
\newcommand{\spic}{\text{Spin$^c$}}

\newcommand{\ov}{\overline}

\newcommand{\bh}{\mathcal H}
\newcommand{\bhd}{\mathcal H_+^2}

\newcommand{\bl}{\mathcal{L}}
\newcommand{\bg}{\mathcal{G}}
\newcommand{\bm}{\mathcal{M}}
\newcommand{\ba}{\mathcal {A}}
\newcommand{\naz}{\nabla_{A_0}}
\newcommand{\na}{\nabla_{A}}

\newcommand{\clif}{\text{End}{(S)}}

\newcommand{\ab}{_{\alpha \beta}}
\newcommand{\gab}{g_{\alpha \beta}}
\newcommand{\tgab}{\tilde g_{\alpha \beta}}

\newcommand{\teal}{\tilde e_\alpha}

\newcommand{\Rr}{\mathbb{R}}
\newcommand{\Hh}{\mathbb{H}}
\newcommand{\Zz}{\mathbb{Z}}
\newcommand{\Cc}{\mathbb{C}}

\newcommand{\id}{\text{Id}}
\newcommand{\ka}{\kappa}
\newcommand{\al}{\alpha}
\newcommand{\ep}{\epsilon}
\newcommand{\si}{\mathfrak{s}}

\newcommand{\om}{\omega}
\newcommand{\Om}{\Omega}

\newcommand{\Ga}{\Gamma}
\newcommand{\Lam}{\Lambda}
\newcommand{\lam}{\lambda}
\newcommand{\kah}{\text{K\"{a}hler}}
\newcommand{\del}{\bar \partial}

\newcommand{\beq}{\begin{eqnarray*}}
\newcommand{\eeq}{\end{eqnarray*}}
\newcommand{\bpr}{\begin{preuve}}
\newcommand{\epr}{\end{preuve}}

\newenvironment{preuve}[1][]
{\vskip 2mm\noindent {\it  D\'emonstration #1.   }}{$\Box$ \vskip 2mm}

\newtheorem{defi}{D\'efinition}
\newtheorem{theo}{Th\'eor\`eme}
\newtheorem{lemme}{Lemme}
\newtheorem{prop}{Proposition}
\newtheorem{coro}{Corollaire}

\title{Invariants de Seiberg-Witten et courbes r\'eelles}
\author{Damien Gayet			}

\begin{document}
%\large
\maketitle

\centerline{\textbf{Abstract}}
On a compact oriented four-manifold with an orientation
preserving involution $c$, we count  solutions of 
Seiberg-Witten equations, which are moreover symmetrical
in relation to $c$, to construct  ``real'' Seiberg-Witten
invariants. Using Taubes' results, we prove that 
on a symplectic almost complex manifold with an antisymplectic
and antiholomorphic involution, this invariants are
not all trivial, and that the canonical bundle
is represented by a real holomorphic curve.

\centerline{\textbf{R\'esum\'e}}
Sur une variété compacte de dimension 4 orientée possédant
une involution $c$ préservant l'orientation, nous comptons
les solutions symétriques par rapport à $c$
des équations de Seiberg-Witten pour définir des
invariants ``réels'' de Seiberg-Witten. Dans le
cas d'une variété symplectique presque complexe munie d'une
involution antisymplectique et antiholomorphe, nous utilisons les résultats
de Taubes pour démontrer d'une part que les invariants 
réels ne sont pas tous nuls, d'autre part que le fibré canonique
de la variété est représenté par une courbe réelle holomorphe.\\

\textsc{Code mati\`ere AMS}: 14P25, 53D05, 57R57.

\thispagestyle{empty}
\subsection{Introduction}
Peu après l'apparition dans [Wi] des invariants de
Seiberg-Witten, Taubes montre que ces invariants
ne sont pas tous triviaux sur une variété symplectique. 
Il en déduit, entre autres, l'existence 
de courbes $J$-holomorphes repr\'esentant 
le fibré canonique de la variété.

Le but de cet article est d'\'etablir 
une version \textit{réelle} de ces résultats, c'est à dire
d'obtenir des courbes non seulement holomorphes, mais
aussi invariantes par une structure réelle sur la variété
symplectique. Plus précisément, considérons une
variété symplectique réelle $(X,\om,J,c)$, où
$\om$ et une forme symplectique, $J$ une structure
presque complexe compatible avec $\om$, et $c$ une
involution $J$-antiolomorphe et antisymplectique. Nous
démontrons le résultat suivant :
\begin{theo}
Soit $(X,\om,J,c)$ une variété de dimension 4 symplectique compacte réelle, 
vérifiant $b_+ >1$. Alors 
le fibré canonique est représenté par une courbe $J$-holomorphe
invariante par $c$, 
a priori singulière.
\end{theo}

Pour cela,  nous nous plaçons plus généralement
dans une variété riemannienne
$(X,g)$ de dimension 4 orientée munie d'une involution $c$ isométrique
et préservant l'orientation. 
Nous définissons d'abord la notion de 
 structure de $\spic$  $c$-réelle. Sur ces structures,
la différentielle $dc$ 
agissant sur le fibré des repères orthonormés orientés
se relève en une involution antilinéaire.
Cette structure est nécessaire pour pouvoir construire
des courbes réelles. En effet, 
dans une
variété presque complexe réelle,
  la structure standard tordue
par  un fibré en cercles 
représenté par  une courbe $J$-holomorphe
réelle  est $c$-réelle. Il est facile
de constater que la structure standard 
tordue par le fibré canonique vérifie cette
condition nécessaire.

Nous consid\'erons alors 
les solutions ``symétriques'' par rapport à $c$ 
des équations de Seiberg-Witten. 
L'espace des solutions,
tout comme dans le cas classique, jouit de propriétés
d'invariance, de lisséité, de compacité, d'orientabilité,
 et permet de définir un invariant
$SW(\si,c)\in \Zz$ indépendant de la métrique pour laquelle
$c$ est isométrique, et vérifiant pour tout
difféomorphisme $f$ de la variété $SW(f^*\si, fcf^{-1}) = SW(\si,c)$.

Dans le cas d'une variété presque complexe symplectique réelle,
la démonstration par Taubes du cas classique s'adapte à notre situation
et permet de démontrer que
 l'invariant réel est non nul
pour la structure standard, twisté ou non
par le fibré canonique.

Pour construire  une courbe $J$-holomorphe représentant
le fibré canonique,
Taubes
utilise  une suite de déformations
des équations classiques. L'invariant étant indifférent
à ces déformations, sa non trivialité implique
à chaque étape l'existence d'une solution. En remarquant
qu'une solution définit entre autres une section 
du fibré canonique (dans le cas twisté), Taubes
montre ensuite que que la suite du lieu des   zéros  de ces 
sections converge
vers une courbe $J$-holomorphe. 
Sans aucun effort supplémentaire que celui
de constater que cette suite de zéros est invariante 
dans notre cas par l'involution $c$, nous obtenons le 
théorème 1. 
\paragraph{Perspectives.}
Taubes [Ta2] a en fait établi que les invariants de Seiberg-Witten 
sur une variété symplectique sont égaux à un certain type 
d'invariants de Gromov-Witten. Ceux-ci compte le nombre
de courbes $J$-holomorphes de degré fixé 
passant par un certain nombre
de points fixes. 
Dans un  travail \`a venir, nous espérons montrer de façon analogue que  
les invariants réels que nous avons construits sont égaux 
 \`a un certain type d'invariants de  Welschinger [We], comptant
des courbes $J$-holomorphes réelles.

\paragraph{Résumé des parties de l'article.}
Dans la première partie, nous construisons les invariants
réels de Seiberg-Witten, dans un cadre général. Nous
expliquons d'abord quelles sont les structures de $\spic$
$c$-réelles. Ensuite nous montrons que l'involution
se relève aux objets naturellement associés à une structure 
de $\spic$ $c$-réelle : espaces de spineurs, espace des connexions
sur l'espace des spineurs. Nous définissons les équations réelles
de Seiberg-Witten, puis nous étudions les propriétés
de l'espace des solutions, en particulier sa dimension.
Dans la deuxième partie, nous montrons d'une part
que dans une variété symplectique réelle, les
invariants construits ne sont pas tous nuls, et d'autre part  que
quand l'invariant d'une structure de $\spic$ est non nul, 
on peut associer à celle-ci une courbe $J$-holomorphe
réelle.

%%%%%%%%%%%%%%%%%%%%%%%%%%%%%%%%%%%%%%%%%%%%%%%%%%%%%%%%%%%%%%%%%%%%%%%%%%%%%%%%%%%%%%%%%%%%%%%%%%%%
		%Construction d'un invariant r\'eel de Seiberg-Witten
%%%%%%%%%%%%%%%%%%%%%%%%%%%%%%%%%%%%%%%%%%%%%%%%%%%%%%%%%%%%%%%%%%%%%%%%%%%%%%%%%%%%%%%%%%%%%%%%

 \section{Construction d'un invariant r\'eel de Seiberg-Witten}

Soit $X$ une vari\'et\'e lisse de dimension 4 orient\'ee,
 munie d'une involution  $c$ lisse pr\'eservant l'orientation.
Il existe toujours  une structure riemanienne 
compatible avec cette situation topologique :
\begin{lemme} L'ensemble $\text{Met} (X,c)$ des 
 m\'etriques riemaniennes  pour lesquelles $c$ est une
isom\'etrie est un convexe non vide. 
\end{lemme}
\begin{preuve} Il est clair que $\text{Met} (X,c)$ est 
un convexe. Montrons qu'il est non vide. L'ensemble des métriques 
riemannienne sur $X$ de volume $1$
est un convexe compact, et   
l'application continue qui à une métrique $g$  associe le tiré en arrière
$c^*g$ laisse  cet ensemble invariant. En effet, 
$$\int_X vol_{c^*g}(x) = \int_X \det d_xc\ . vol_{g}(c(x))= 
\int_X \det d_x c\ .  \det d_{c(x)}c\ . vol_g(x) = \int_X vol_g(x).$$
Le théorème de Schauder implique maintenant le résultat.  
\end{preuve}

		%%%%%%%%%%%%%%%%%%%%%Structures Reelles%%%%%%%%%%%%%%%%%%%%%

\subsection{Structures réelles}
\subsubsection {Une remarque éclairante pour la suite}
Soit $f: \Cc^2 \to \Cc$ une application lisse définie
sur le plan complexe. Si $f$ est holomorphe,
l'involution antilinéaire 
$\kappa (f)$ définie par $\ka( f) (z) = \bar f (\bar z)$ 
est également
holomorphe. Si de plus $f$ est invariante
par $\ka$, alors le 
lieu des zéros de $f$  est une sous-varété analytique complexe réelle,
c'est à dire invariante par conjugaison.
La suite de cette section explique comment construire
des équivalents de $\ka$ sur 
un fibré en droites, puis
sur les  objets associés à une  structure de $\spic$. 
%***********Structures de Spin r\'eelles*************************

%
%
%*Structures réelles et fibré en droites complexes********************
\subsubsection {Fibrés en cercles  $c$-réels}
Soit $Q$ un fibré principal en cercles.
Il nous faut savoir 
quand il est possible de 
 relever $c$ en une involution antilinéaire sur $Q$.
Pour cela, nous avons besoin   des groupes et du faisceau suivants :
\beq
\mathcal G & = &  \{f \in \Ga(X,\un)\}\\
r\mathcal G & = &  \{f \in \Ga(X,\un), \ f\circ c = f\}\\
r\mathcal G_\al &=& \{f \in \Ga(U_\al,\un), \ f\circ c = f\}
\eeq
\textbf{Remarque} : Sans perte de généralité, on a supposé et on supposera
 toujours que les ouverts
de trivialisation $U_\al$  sont invariants par $c$.
% Nous voulons donc construire une involution
%$$\kappa_L: \Ga(X,L) \rightarrow \Ga(X,L)$$ vérifiant pour tous $\psi \in \Ga(X,L)$
%et $f \in \Ga (X, \un)$
 
%$$\kappa_L(f\otimes \psi)(x) = \bar f (c(x)) \otimes c_L(\psi) (x).$$ 
%Si $L$ admet une telle structure, on dira que $L$ est $c-$réel.
%Pour cela,  nous définissons
%le sous-groupe $$rG = \{f \in \Ga(X,\un), \ f\circ c = \bar f\}.$$
	\begin{lemme}
Soit $Q$ un \un-fibré principal sur $X$. 
L'application $c:X \to X$ se relève en une application
$c_Q: Q \to c^*Q$ vérifiant $$c_Q(\lam \psi) = \bar \lam c_Q(\psi)$$
pour toute application  $\lambda \in  \mathcal G$ et
toute section $\psi \in \Ga(Q)$
si et seulement si $c^*Q\otimes Q$ est trivial en tant qu'élément
de $H^1(X,\bg)$. Il existe une telle application involutive
si et seulement si $c^*Q\otimes Q$ est trivial en tant qu'élément
de $H^1(X,r\bg)$. Dans ce dernier cas,
l'ensemble  de telles involutions est  isomorphe
à $r\mathcal G$.
	\end{lemme}
\begin{preuve}
 Si $h\ab \in \un$ sont les fonctions
de transition de $Q \in H^1(X,\un)$, le fibré
$c^*Q\otimes Q $ a pour fonctions de transitions
$ h\ab(c)h\ab$. Ces fonctions sont des éléments du faisceau $r\mathcal G$,
et donc définissent un élément de $H^1(X,r\bg)$. 

Maintenant, soient $\ep_\al$ des trivialisations locales
de $Q$, telles que 
$\ep_\al = h\ab \ep_\beta.$ 
Si l'involution se relève,  il existe des fonctions
 $\ka_\al \in \Ga (U_\al, \un)$
telles que 
$ c_Q( \ep_\al) = \ka_\al \ep_\al(c)$. 
On a  alors $$\bar h\ab (c) \ka_\beta = \ka_\al h\ab,$$
 ce qui
montre que le cocycle
$\{h\ab(c)h\ab\}$ est nul dans $H^1(X,\un)$, 
et donc que $c^*Q = -Q$. Si maintenant
l'application $c_Q$ est une involution,
on a $\ka_\al \in r\mathcal G_\al$, ce qui signifie cette 
fois que le cocycle
$\{h\ab(c)h\ab\}$ est nul dans $H^1(X,r\bg)$.
 La réciproque suit la voie inverse. Enfin, soit $c_Q$ un relèvement
involutif de $c$ à $Q$. Il est facile de vérifier que toute autre relèvement
involutif est un produit  de $c_Q$ par un élément de $r\mathcal G$.
\end{preuve}

\begin{defi} Un fibré $Q$ est dit $c$-réel s'il existe une
involution antilinéaire de $Q$ dans $c^*Q$.
\end{defi}

 Les propositions suivantes montrent que cette définition est naturelle :
\begin{prop} Soit $(X,J)$ une variété presque complexe, 
et $c$ une involution $J$-antiholomorphe. 
Si $C$ est une courbe $J$-complexe plongée 
dans $(X,J)$, invariante par $c$, alors  le \un-fibré associé est $c$-réel. 
\end{prop}
\begin{preuve}
Supposons d'abord que $J$ est intégrable. Soient
$(f_\al)$ des applications holomorphes irr\'eductibles
définissant localement
$C$. L'antiholomorphie de $c$ implique que $\bar f_\al (c)$
est holomorphe, et l'invariance
de $C$ par $c$  implique que $f_\al$ et 
$\bar f_\al (c)$ s'annulent en m\^eme temps. Il existe donc des applications
$g_\al$ holomorphes ne s'annulant pas et des entiers positifs $m$,
telles que $\bar f_\al (c) = g_\al f_\al^m$. 
On a donc aussi 
$$f_\al = \bar g_\al(c) \bar f_\al^m(c)=  \bar g_\al(c)g_\al^m f_\al^{m^2}.$$
L'application  $\bar g_\al(c)g_\al^m$ est holomorphe,
donc l'\'egalit\'e  n'est possible que si $m=1$. De plus 
$g_\al \bar g_\al(c)$ vaut 1
en dehors de $C$ et donc partout. On a donc $g_\al/|g_\al| \in r\bg_\al$.
Maintenant,
le cocycle $$\{ h\ab = \frac{f_\al}{|f_\al|}\frac{|f_\beta|}{f_\beta} \}$$
définit le fibré en cercles 
  associé à  $C$, et puisque $h\ab(c)h\ab = g_\beta/|g_\beta|. |g_\al|/g_\al$, 
le fibré est $c$-réel. 

Le cas général se ramène au cas intégrable après
avoir remarqué qu'il est possible de construire une
structure presque complexe intégrable sur $U_\al$ 
égale à $J$ aux points de $C$.
\end{preuve}
\begin{prop}Soit $(X,J)$ une variété presque complexe munie
d'une involution $J$-antiholomorphe. Alors le fibré anticanonique
$K^{-1}$ est $c$-réel.
\end{prop}
\bpr
Le fibré $K^{-1}$ est identifié aux 
2-covecteurs complexes de type $(0,2)$. 
L'application $c_{K^{-1}}$ qui envoie
$\om$ sur $\overline{\om} \circ dc^{-1}$ 
relève $c$, préserve $K^{-1}$, est antilinéaire et involutive.
\epr  
\subsubsection {Une involution sur $\spic(X)$ }
Un fibré en cercles $Q$ est $c$-réel
si, entre autres, $c^*Q = -Q$.
 Avant de définir la notion de structure de $\spic$ $c$-réelle,
nous rappelons la construction de l'involution sur $\spic(X)$
envoyant une structure de $\spic$ $\si$ de fibré déterminant 
$Q$ sur la structure de $\spic$ $-\si$ de fibré déterminant $-Q$.

\paragraph{Discussion locale.}
Ce paragraphe reprend l'expos\'e de Morgan dans son livre 
[Mo], pp. 100-102.
D\'efinissons d'abord la conjugaison complexe 
sur  l'alg\`ebre de Clifford complexe
$Cl(4)\otimes \Cc$. Elle applique
simplement $e$ sur $\bar e$. Par induction sur 
$\spic(4)\subset Cl(4)\otimes \Cc$, on obtient 
la conjugaison sur $\spic(4)$.

On sait que  $ Cl(4)$ est isomorphe en tant qu'alg\`ebre
\`a $\text{Mat}(2,\Hh)$, o\`u $\Hh$ est le corps des quaternions.
Cet ensemble agit sur $\Hh\oplus \Hh$ par multiplication 
matricielle quaternionique \`a gauche, tandis que l'ensemble
des complexes $\Cc$ agit \`a droite par multiplication \`a droite par 
$i$. Maintenant, la multiplication \`a droite par $j$ 
d\'efinit une application $\tau$ de $\Hh\oplus \Hh$ dans lui-m\^eme
v\'erifiant pour tout $e\in Cl(4)\otimes \Cc$ et tout $s \in \Hh\oplus \Hh$ 
$$ 
\tau (e\cdot s) = \bar e \cdot \tau (s).
$$

\paragraph{Globalisation.}
Soit maintenant $P$ le $\SO$-fibr\'e principal des rep\`eres  orthonorm\'es
directs. 
Si $\si$ une structure de $\spic$, $\si$ est un rev\^etement double
de $P\times_X Q$, où $Q$ est le $\un$-fibré déterminant de $\si$. 
L'application induisant l'identité sur $P$ et la conjugaison
complexe sur $Q$ définit par tiré en arrière 
une nouvelle structure de $\spic$ notée $-\si$. On remarque
que le fibré déterminant de $-\si$ est $-Q$.
Il existe alors une application antilin\'eaire 
$$ \tau : \si \to -\si.$$
Cette application  d\'efinit un isomorphisme
antilin\'eaire
$$
\tau_S: S(\si) \to S(-\si),
$$
o\`u $S(\si)$ est le fibr\'e de spineurs 
$\si \times_{\spic(4)}\Hh\oplus \Hh = S^+ \oplus S^-$. On a de plus
par construction
$$\tau_S (e\cdot\psi) = \bar e \cdot  \tau_S(\psi)$$
pour toute section $e$ de $Cl(TX)\otimes \Cc\simeq TX\otimes \Cc$ et
toute section $\psi$ de $S(\si)$. Enfin, $\tau_S$ préserve
la décompostion en spineurs positifs et négatifs.
%*Structures de $\spic$ r\'eelle**********************
\subsubsection {Structures de $\spic$ r\'eelles}
La différentielle   $dc$ d\'efinit un morphisme
$dc: P \to  c^*P$ relevant $c: X \to X$. Appelons 
$\mathcal C_\si$ l'ensemble des rel\`evements de $dc$
\`a $\si$ antilin\'eaires.
	\begin{prop}
Soit $\si$ une structure de $\spic$ et $Q$
son fibr\'e d\'eterminant. Si $Q$ est $c$-réel,
$\si$ d\'efinit naturellement un \'el\'ement $\mu$ de $H^1(X,\Zz_2)$, tel que
 $\mu = 0$ si et seulement si $c^*\si = -\si$. Dans ce cas 
 $\mathcal C_\si \not= \emptyset$, et il existe 
un élément $\ep(\si) \in \Zz_2$ tel que $\forall c_\si \in  C_\si$,
$c_\si ^2 = \ep(\si) \id$. 
	\end{prop}
	\begin{preuve}
Soit $e_\al$ une trivialisation de $P$,
et $\gab$  les fonctions de transition
 à valeurs dans $\SO$ vérifiant $e_\beta = e_\al \gab$.
Nommons de plus  $\teal$ des trivialisations
du fibré $\si$, 
 $\tgab$ les fonctions
de transitions de $\si$ et $h\ab = \det \tgab$ celles de $Q$.
Supposons enfin que les $\tilde e_\al$(resp. $\tgab$) sont des relev\'es de
$e_\al$ (resp. $\gab$).
Soit
 $k_\al$ des applications à valeurs dans  $\SO$, telles que
$$ d_xc (e_\al) =  e_\al (c(x))k_\al^{-1}(x),$$
et $\ka_\al$ des applications appartenant à  $r\mathcal G_\al$
telles que 
$$
h\ab (c) = \ka_\al^{-1}\bar h\ab \ka_\beta.
$$
Choisissons des relev\'es arbitraires $\tilde k_\al$
\`a $\spic(4)$ des fonctions $k_\al \times \ka_\al$.
Les  cocycles
$$
\tilde k_\al^{-1}(\tgab \bar h\ab) \tilde k_\beta
$$
et $\tgab(c)$ ont m\^eme image par projection
sur $\SO\times \un$, si bien qu'ils ne diffèrent
que d'un signe $\mu\ab$.  Leur rapport définit donc
un élément $\{\mu\ab\} \in H^1(X,\Zz/2\Zz)$.
Si $\mu$ est nul,
il existe des constantes $\eta_\al \in \{\pm 1\}$
telles que $\mu\ab = \eta_\al/\eta_\beta$. Alors 
l'application  
$$c_\si ( \tilde e_\al) =  \tilde e_\al(c) \tilde k_\al^{-1}\eta_\al$$
est bien d\'efinie.
En projetant sur les deux facteurs $\SO$ et $\un$, on constate
que $c_\si^2 = \pm \id$. 
	\epr
		\begin{defi}Sous les hypothèse de la proposition précédente, 
on dira que $\si$ est $c$-réelle si $\ep(\si) = 1$.
		\end{defi}

%Structure réelle et fibré de spineurs*********************
\subsubsection{Structure réelle et fibré de spineurs}
		\begin{prop} Soit $\si$ une structure
de $\spic$ $c$-réelle. Alors il existe une isom\'etrie antilinéaire
$\kappa_S: \Gamma (S) \rightarrow \Gamma (S)$ 
stable pour les deux sous-fibrés $S^+$ et $S^-$,
et induisant sur $\Ga(\clif)\simeq \Om^*(X,\Cc)$ le tiré en arrière barré :
$$\kappa_\Om (\omega) = \ov{c^*\omega}.$$
		\end{prop}
	\begin {preuve} 
On a vu qu'il existait une application anticomplexe 
$\tau_S: S(\si) \to S(-\si)$ induisant
la conjugaison sur les formes complexes sur $X$.
Puisque $\si$ est réelle,
on a d'autre part $-\si = c^*\si$, ce qui donne
un isomorphisme 
$\nu_S: S(-\si) \to S(c^*\si)$ tel que pour tout champ
de vecteur $e$,
$$
\nu_S (e\cdot\psi) = dc(e)\cdot\nu_S(\psi).
$$
Au total, on obtient une application 
$c_S: S(\si) \to S(c^*\si)$
telle que l'application 
\beq
\ka_S: \Gamma (S) &\to& \Gamma (S) \\
\ka_S (\psi) &=& c_S \circ \psi \circ c
\eeq
v\'erifie, en identifiant de facon naturelle $TX$ avec
$T^*X$, 
$$ \ka_S (\om \cdot \psi)= \ka_\Om(\om)\cdot \ka_S(\psi).$$
 Enfin par construction l'application laisse stable
les deux composantes $S^\pm$ de $S$. 
		\end{preuve}

%*Involution sur l'espace des connexions**************************
\subsubsection{Involution sur l'espace des connexions}
Commençons  par l'espace $\ba(Q)$ 
des connexions unitaires
d'un $\un$-fibré $Q$.
		\begin{prop}Soit $Q$ un $\un$-fibré principal $c$-réel.
		Alors il existe une involution antilinéaire
		$\ka_\ba$ de $\ba(Q)$ dans lui-m\^eme, telle 
		que la courbure $F_{\ka_\ba(A)}$ vérifie
		$$F_{\ka_\ba(A)}= \ka_\Om (F_A).$$
		\end{prop}
	\begin{preuve}
Soit$A$ une connexion hermitienne sur $Q$. 
 Si $\ep_\al$ est une trivialisation
locale de $Q$, il existe des 1-formes
$A_\al$ à valeurs dans $\Rr$
telles que $A \ep_\al = iA_\al \ep_\al$.
Par ailleurs, puisque $Q$ est $c$-réel,
il existe des fonctions $\ka_\al$ à valeurs 
dans $r\mathcal G_\al$ telles que $c_Q \ep_\al (c) =  \ep_\al \kappa_\al^{-1}$.
Soit $  \ka_\ba(A)$ la connexion définie par 
$$ \ka_\ba(A) \ep_\al = (-i c^*A_\al + d \log \ka_\al)\ep_\al.$$
On vérifie facilement que cette connexion est bien définie
et qu'elle est involutive. 
 La courbure de $A$
est $F_A = d(iA_\al) =d(iA_\beta)$, et donc la courbure de son image
est $d(-ic^*A_\al) = \ka_\Om (F_A)$.
	\end{preuve}

Rappelons la proposition suivante :
		\begin{prop} Soit $\si$ une structure de $\spic$ et $Q$ son
		fibré déterminant, ainsi qu'une connexion unitaire $A$ sur
		$Q$. Alors il existe une unique connexion $\na$
		sur l'espace des spineurs
		 induisant la connexion de Levi-Cività $\nabla$  sur $P$ et 
		la connexion $A$ sur $Q$.
		\end{prop}

%Supposons maintenant que $\si$ soit $c$-réelle. La proposition
%précédente permet d'identifier l'espace des 
%connexions levicivitiennes sur l'espace des spineurs
%$\mathcal A (\si)$avec 
%$\mathcal{A}(Q)$, où $Q$ est le fibré déterminant de $\si$.
% La proposition
%précédente nous permet de définir 
%		\begin{defi}
%	Soit $\si$ une structure de $\spic$ réelle, et $Q$ son
%	fibré déterminant. On définit $\kappa_\ba : \ba(\si) \to \ba(\si)$
%	de la façon suivante: la connexion $\ka_\ba (\na)$
%	est la connexion levicivitienne induisant la connexion 
%	${\ka_\ba(A)}$ sur $Q$.
%		\end{defi}
		\begin{lemme}
		Soit $\si$ une structure de $\spic$ $c$-r\'eelle, 
		et $Q$ son fibr\'e d\'eterminant.
		Pour toute connexion $A \in \ba(Q)$
		on a les relations de commutation suivantes :
		\beq
 		\nabla = \ka_\Om \nabla \ka_\Om,\\
		\nabla_{\ka_\ba(A)} = \ka_S \na \ka_S.
		\eeq
		\end{lemme}
	\begin{preuve}
Soient $(e_i)$ une section orthonormée de $TX$ sur $U_\al$.
La connexion de Levi-Cività s'écrit 
$$\nabla e_i = \sum_{i=1}^{4} \om_{ij}e_j,$$
où $\om_{ij}$ sont des 1-formes à valeurs réelles.
Composons par le tiré en arrière de $c$ :
$$c^* \nabla e_i = \sum c^*\om_{ij} c^*e_j. $$
Nous avons besoin du lemme suivant :
		\begin{lemme}
		Soient $e'_i =c^*e_i$.
		 Si $\om_{ij}'$  est l'ensemble
		des 1-formes définissant la connexion de Levi-Cività
		dans la base orthonormée $e'_i$, 
		alors on a $\om'_{ij} = c^*\om_{ij}$.
		\end{lemme}
	\begin{preuve}En effet, si $\nabla e_j (e_i)= \sum \Ga_{ij}^{k}e_k$,
avec 
$$\Ga_{i,j}^{k} = \frac{1}{2}g^{il}(g_{jl,k}+ g_{kl,j} - g_{jk,l}),$$
où $g^{il}$ définit la métrique duale et $g_{jl,k} = d g_{jl}(e_k)$.
Dans la trivialisation $e'_i$, les $g'_{ij}$ égaux
à $g_{ij}$ car $c$ est une isométrie, et donc $g'_{jl,k} = c^*g_{jl,k}$.
Au total, on constate que $\Ga_{ij}^{k'} = c^*\Ga_{ij}^k.$
	\end{preuve}
Gr\^ace à ce lemme, on obtient donc aisément que
$\ka_\Om \nabla \ka_\Om = \nabla.$
Maintenant soit 
$A\in \ba(Q)$. La connexion
spinorielle  $\na$ s'écrit en coordonnées
$$ \na = A_\al \id + \sum_{i=1}^{4}\om_{i,j} e_i \cdot e_j,$$
où $\cdot$ est la multiplication de Clifford.
Gr\^ace au résultat sur la connexion de Levi-Cività,
il est alors aisé de constater que $	\nabla_{\ka_\ba(A)} = \ka_S \na \ka_S.$
	\end{preuve}

%%%%%%%%%%%%%%%%%%%%%%Equations réelles de Seiberg-Witten%%%%%%%%%%%%%%%%%%%%%

\subsection{Equations réelles de Seiberg-Witten}
Les équations que nous définissons sont celles de Seiberg-Witten,
mais nous imposons aux inconnues, c'est-à-dire 
un spineur positif et une connexion, d'\^etre invariants
par l'involution réelle. Rappelons, avant de décrire
l'espace des solutions, les équations classiques.

%Les équations de Seiberg-Witten**********************

\subsubsection{Les structures supplémentaires}
\paragraph{Formes autoduales.}
Rappelons que si $\{e_1, e_2, e_3, e_4\}$ est une
base orthonormale locale directe, l'étoile de Hodge est définie
de la façon suivante :
\begin{eqnarray*}
\ast : \Lam^2 TX \to \Lam^2 TX\\
e_i\land e_j \to e_r \land e_s
\end{eqnarray*}
où  $(i,j,r,s)$ est une permutation paire de (1,2,3,4).
Par ailleurs  
via l'isomorphisme entre  $\Lam^*TX\otimes \Cc$ et $\clif$,
le sous-espace $\Lam_+^2TX\otimes i\Rr$  s'identifie avec
 les endomorphismes de $S^+$
hermitiens de trace nulle.

\paragraph{Opérateur de Dirac.}
Soit $\psi$ un spineur, c'est à dire une section
du fibré $S$, et $\na$ une connexion
sur $S$. 
La dérivée $\na \psi$ est une section  de $\Ga( T^*X \otimes S)$.
On peut, en utilisant l'identification entre $T^*X\otimes \Cc$ et $\clif$, 
définir maintenant l'opérateur de Dirac 
$\di: \Ga( S) \to \Ga( S).$
L'opérateur envoie une section de $S^\pm$
sur une section de $S^\mp$. 
En coordonnées, on trouve 
$$ \di \psi = \sum_{i=1}^{4} e_i\cdot \nabla_{e_i}\psi.$$

\paragraph{L'application quadratique $q$.}
Soit $\psi$ un spineur positif, c'est à dire 
une section de $S^+$. 
Le spineur définit une 2-forme autoduale $q(\psi)$ \`a valeurs
complexes pures (ou
un endomorphisme de $S^+$ hermitien de trace nulle) en posant
$$ q (\psi) = \psi \otimes \bar \psi - \frac{|\psi|^2}{2} \id.$$
Dans une base locale de $S^+$, si $\psi = (\al,\beta)$, on a 
$$ q (\psi) = \frac{1}{2}\left( \begin{array}{cc}
				|\al|^2-|\beta|^2  & 2\alpha \bar \beta\\
				2\bar \alpha \beta   &-( |\alpha|^2-|\beta|^2)
				\end{array} \right).
$$
\subsection{Les équations}
Soit $X$ une variété compacte orientée riemannienne munie
d'une structure de $\spic$. On appelle $Q$ 
le fibré déterminant associé, et $\Ga^+\times \ba$ le produit
$\Ga( S^+) \times \ba(Q)$.  Pour toute 2-forme $h$ autoduale, 
les équations  perturbées  de Seiberg-Witten
ont pour inconnues un spineur positif  $\psi$ 
et une connexion unitaire  $A$ sur $Q$ :
\begin{eqnarray*}
(\psi, A) & \in & \Ga^+\times \ba\\
\di \psi& =& 0\\
F^+_A &=& q(\psi)+ ih,
\end{eqnarray*}
où $F^+_A$ est la partie autoduale de la courbure de 
la connexion $A$. 
\vskip 2mm

\subsection{Les équations réelles}
Les équations réelles de Seiberg-Witten que
nous définissons sont les m\^emes, sauf que l'on 
impose à $\psi$ et $A$ d'\^etre
invariants par les involution $\ka_S$ et
$\ka_\ba$ construites précédemment; 
si 
\beq
r \Ga^\pm &= &\{\psi \in \Ga^\pm, \ \ka_S \psi = \psi\}\\
r \ba & = & \{A \in \ba(Q), \ \ka_\ba(A) = A\}\\
a \Om^* & =& \{\om \in \Om^*(X,\Rr),\  \ka_\Om \om = - \om\}
\eeq
les équations  $E_h(\si, c)$ sont définies par  
\begin{eqnarray*}
(\psi, A) &\in &r\Ga^+\times r\ba\\
\di \psi& = &0\\
F^+_A &= &q (\psi)+ ih,
\end{eqnarray*}
où $h$ est une 2-forme auto-duale anti-$c$-réelle, 
c'est à dire $h \in a \Om^2_+  $. 
On notera $M_h (\si,c)$ l'ensemble
des solutions de $E_h(\si,c)$.
\paragraph{Groupe de jauge.}
Le sous-groupe $r\mathcal G \subset \mathcal G$ 
agit sur l'espace $r\Ga^+\times r\ba$
de la façon suivante :
$$ f \ .\  (\psi, A) = (g\psi,\   A +gd(g^{-1})).$$
Il est facile de vérifier qu'il 
laisse invariant  $M_h(\si, c)$. Par ailleurs
le stabilisateur d'un élément $(\psi, A)$ est trivial
si $\psi$ n'est pas le spineur nul, et est égal à $\Zz_2$
sinon. Dans le premier cas on dira que le couple
$(\psi, A)$ est irr\'eductible.
 Enfin, on note 
$\bm_h (\si,c)$ l'ensemble des classes d'équivalences
de $M_h(\si,c)$ sous l'action de $r\mathcal G$.

%%%%%%%%%%%%%%%%%%%%%%L'espace de modules%%%%%%%%%%%%%%%%%%%%%

\subsection{Propriétés de l'ensemble des solutions}

\subsubsection{Le théorème principal}
%Avant d'\'enoncer le th\'eor\`eme principal, nous devons d\'efinir 
%des les équivalents de $b^1$ et $b^+$ pour notre situation:

	\begin{defi}
	On appelle $a\bh^1$ (resp. $a\bhd$) l'espace
	des 1-formes (resp. 2-formes autoduales) 
	réelles harmoniques anti-$c$-réelles, i.e 
	$a\bh^1= \bh^1 \cap a\Om^1$ (resp. $a\bhd = \bhd \cap a\Om^2$).
	\end{defi}

		\begin{theo}
		Soit $(X,\ c,\ g)$ une variété réelle compacte 
		telle que $\dim a\bhd >0$. 
		Alors pour toute 2-forme 
		$h$ anti-$c$-réelle autoduale générique,
		et toute  structure 
		de $\spic$ $c$-réelle $\si$,  
		l'espace $\bm_h(\si,c)$ est une variété compacte orientée
		de dimension 
		$$d= \dim a\bh^1 - \dim a\bhd  + \frac{1}{8}(c_1(Q)^2 - \tau),$$
		où $\tau$ est la signature de $X$ et $Q$ le fibré déterminant
		de $\si$.  
		\end{theo}
Si $d$
est strictement positive, choisissons un point base
dans $X$, et soit $r\mathcal G_0$ l'ensemble des 
éléments de $r\mathcal G$ valant $1$ en ce point. 
Soit $\bm^0_h(\si, c)$ le quotient 
de $M_h(\si,c)$ par ce sous-groupe. Alors la projection
$\bm^0$ sur $\bm$ définit un fibré $\un$-principal. 
Notons $\mu$ la première classe de Chern de ce fibré.
Nous pouvons maintenant définir les invariants
réels de Seiberg-Witten. 
	\begin{defi}
	$\bullet$ Si la dimension $d$ ci-dessus est strictement négative ou impaire, 
	$SW_h(\si,c) = 0$.\\
$	\bullet$ Si $d$ est nul, $SW_h(\si,c)$ 
	est le nombre 
	le nombre d'éléments de $\bm_h(\si, c)$ comptés avec 
	leur signe.\\
	$\bullet$ 
	Si $d$ est strictement positive et paire,  	
	$$SW_h(\si,c) = \int_{\bm_h(\si,c)}\mu^{d/2},$$
	\end{defi}
		\begin{prop}
		Si $\dim a\bh^2_+>1$, l'invariant défini ci-dessus 
		est indépendant de la perturbation $h$ 
		et de la métrique pour laquelle $c$ est isométrique.
		On appelle $SW(\si,c)$ cet entier relatif.
		De plus, si $f$ est un difféomorphisme sur $X$ on la relation :
		$$SW(\si, fcf^{-1}) = SW(f^*\si, c).$$
		\end{prop}

Nous allons démontrer ce théorème en plusieurs temps. 
D'abord nous analyserons les problèmes d'indices, 
ensuite celui de la transversalité, et enfin celui 
de la compacité.

\subsubsection{L'indice}
Le lemme suivant, basé sur les relations
de covariances du lemme 3,
 établit la covariance des équations
de Seiberg-Witten sous l'action des involutions $\ka_S$ et$\ka_\ba$.
		\begin{lemme}
		Soit l'application $f$ définie par 
		 \beq
		f : \Ga^+\times \ba \times \Om^2_+ &\to&\Ga^- \times i\Om^2_+\\ 
		(\psi,\ A,\  h)&\mapsto& (\di \psi,\  F^+_A -q (\psi)-i h).
 		\eeq
		On a alors $ f \circ (\ka_S,\ka_\ba,\ka_\Om) = (\ka_S, \ka_\Om)\circ f$. 
		On peut donc  définir par restriction
			l'application :
		$$ f_r :
		r\Ga^+\times r\ba \times a\Om^2_+ \to r\Ga^- \times ia\Om^2_+ ,$$
		\end{lemme}
	\bpr[du lemme]Rappelons en effet que 
$\na \psi$ est à valeurs dans $T^*X \otimes S^+$, et
$\di$ fait simplement  agir 
multiplication de Clifford  de $T^*X$ sur $S^+$. On a donc
$$\di = A_\al \cdot\id + \sum_{i=1}^{4}\om_{i,j}\cdot e_i \cdot e_j.$$
Gr\^ace au lemme 3, on obtient  que $\ka_S \di = \di \ka_S$.
Pour la deuxième coordonnée de l'application, il suffit d'abord 
 de remarquer que 
d'une part $F_A =d A_\al$, ensuite que $(\ka_\Om F)^+ =\ka_\Om( F^+)$
puisque $c$ est une isométrie. Montrons
enfin que $q\circ \ka_S = \ka_\Om \circ q$. On a 
$ q (\psi) = \psi \otimes \bar \psi - |\psi|^2 \id$, 
ce qui donne 
$$
q (\ka_S \psi) = \ka_S \psi \otimes \ov{\ka_S \psi} - |\ka_S \psi|^2 \id 
	     = \ka_\Om ( \psi \otimes \bar \psi) - c^*|\psi|^2\ka_\Om (\id)
	     =  \ka_\Om (q (\psi)).
$$	    
\epr
Etudions maintenant la linéarisation de $f_r$. Remarquons d'abord
les égalités suivantes concernant les  tangents :
\beq
%T(\Ga^+\times \ba \times \Om^2_+(X,\Rr))  &=& \Ga^+ \times i\Om^1(X,\Rr) \times \Om^2_+\\
T(r\Ga^+\times r\ba\times a\Om^2_+) &=& r\Ga^+ \times ia\Om^1\times a\Om^2_+\\
T(r \Ga^- \times i a \Om^2_+) & = & r \Ga^- \times i a \Om^2_+.
\eeq
La linéarisation de $f_r$ donne :
\beq
d f_r : T(r\Ga^+\times r\ba\times a\Om^2_+) &\to& T(r \Ga^- \times i a \Om^2_+)  \\ 
        (\psi', a, h')& \mapsto & 
(\di (\psi') -ia\cdot\psi\ ; \ (ida)^+ - 2q(\psi, \psi')-i h').
\eeq
L'algèbre de Lie du groupe de jauge $r\mathcal G$ est par ailleurs
$ia\Om^0$.
Nous avons donc à calculer l'indice du complexe suivant, tiré
de la linéarisation  en une solution
 $(\psi, A)$ de l'action de $r\mathcal G$ et des équations réelles
de Seiberg-Witten :
$$ 
ia\Om^0 \to  r\Ga^+ \times ia\Om^1\times a\Om^2
		\stackrel{df_r}{\longrightarrow} r \Ga^- \times i a \Om^2_+,$$
la première flèche étant donnée par 
$ f \mapsto (-f\psi, \ 2df,\  0).$
L'indice d'un opérateur 
ne dépendant que de la partie principale de l'opérateur,
 le complexe elliptique précédent se résume en la somme
de deux complexes, 
$a\Om^0 \stackrel{d}	{\to} 
a \Om^1  \stackrel{d^+}	{\to} 
a \Om_+^2,$
et
$ 0 \to r\Ga^+ \stackrel{\di}{\to} r\Ga^-.$

Afin de calculer l'indice du second complexe,  on remarque
que la mutliplication complexe à droite sur $\Ga^+$
envoyant $\psi$ sur $\psi.i$ envoie
$r\Ga^{\pm}$ sur $a\Ga^{\pm}$. Cette application
est linéaire et inversible, et commute 
avec l'opérateur $\di$. On obtient donc 
$ \ind D_{A|r\Ga^+}= \frac{1}{2} \ind \di.$

En ce qui concerne l'indice du premier complexe, nous 
avons le 
		\begin{lemme}L'indice du complexe 
		$a\Om^0 \stackrel{d}	{\to} 
				a \Om^1  \stackrel{d^+}	{\to} 
		a \Om_+^2$
		est égal à 
		$ -\dim a\bh^1 + \dim a\bhd$.
		\end{lemme}
	\bpr
L'indice du complexe est égal à la somme alternée 
des dimensions des groupes de cohomologie 
associés à la suite. Le premier groupe est 
simple à calculer : si $f$ est une fonction 
vérifiant à la fois $df = 0$ et $c^*f=- f$, 
$f$ est la constante nulle. 

Si $\al$ appartient à $\coker d_{|a\Om^0} \cap \ker d^+$, alors 
$\al$ est orthogonale à $d(a\Om^0)$. Si $f$ est une fonction
quelconque à valeurs réelles, $f=f_r + f_a$, où
$c^*f_r = f_r$ et $c^*f_a =-f_a$. On a alors
$$(\al, df) = (\al, df_r) = (c^*\al, c^*df_r) 
= (-\al, df_r) =0.$$
Cela signifie que $\delta \al = \star d \star \al = 0$.
Puisque $2 d^+\al = d\al +  d\star \al$, 
 $\al$ est en fait harmonique. Le second groupe de
cohomologie est donc $a\bh^1$.

Calculons maintenant le troisième groupe de 
cohomologie du complexe. Si $\om \in a\Om_+^2$ 
est dans $\coker d^+_{|a\Om^1}$, on démontre
comme précédemment qu'il est en fait
orthogonal à $d^+\Om^1(X,\Rr)$, ce qui
 ajouté à l'autodualité de $\om$ implique
que $\om$ est harmonique, et que le troisième
groupe est $a\bh^2_+$.
	\epr
\subsection{Transversalité}
Tout comme dans le cas classique, on démontre
d'abord gr\^ace au théorème de Sard-Smale que pour 
une perturbation $h$ générique,  $\bm_h(\si,c)$
est une variété lisse de la dimension $d$, aux points \textit{irréductibles}.
L'élément à vérifier dans notre cas est que $df_r$ est surjective
aux solutions  irréductibles de $E_h(\si,c)$.
Ensuite, l'élimination des solutions irréductibles
se fait par le lemme suivant, où la condition
$\dim a \bh^2_+>0 $
est utilisée.

\begin{lemme}
Si $\dim a \bh^2_+ >0 $, 
alors pour toute perturbation $h$ générique,
pour toute structure $\si$ de $\spic$ réelle, 
il n'y a pas de solution réductible de $E_h(\si,c)$.
\end{lemme}
\bpr
S'il existait une solution réductible, on 
aurait $F_A^+ = ih$, avec $A$ et $ih$ 
$c$-réels. En écrivant $A = A_0 + ia$, avec $a\in a\Om^1$
et $A_0$ une connexion base sur le fibré déterminant, on constate
que  l'image par $F_A^+$ de l'ensemble des solutions
réelles réducibles est un espace affine. Si la perturbation
autoduale $ih$ est orthogonale
à cette image, on démontre comme pr\'ec\'edemment que $h$ est 
 harmonique.
La codimension de l'image est donc $\dim a\bh_+^2$.
\epr

\subsection{Orientation}
Nous renvoyons au livre de Morgan [Mo] pour 
la discussion sur l'orientation dans le cas classique,
et qui s'applique \textit{mutatis mutandis} à notre situation. 
Nous ne donnons que la 
\begin{prop}
Le choix d' orientations sur  $a\bh^1$ et
sur $a\bh^2_+$ détermine une orientation
sur $\bm_h(\si,c)$ pour toute perturbation $h$.
\end{prop}

\subsection{Compacité}
La compacité des variétés $\bm_h(\si,c)$ 
se déduit de la compacité du cas classique
en remarquant que la condition d'\^etre
$c$-réel est une condition fermée.

\subsection{Une involution dans la théorie}
Tout comme dans le cas classique, 
on a gr\^ace à l'involution sur $\spic(X)$ décrite au
paragraphe 1.1.3 la 
\begin{prop}
Soit $\si$ une structure de $\spic$ $c$-réelle. Alors
$-\si$ est aussi $c$-réelle, et de plus
$SW(\si, c) = \pm SW(-\si, c)$.
\end{prop}
\bpr
Si $c_\si$ est une involution antilinéaire définie sur $\si$, et
$\tau$ est l'application antilinéaire définie précédemment
de $\si$ dans $-\si$, alors $c_{-\si}= \tau \circ c_\si \circ \tau$
définit une structure réelle sur $-\si$. De plus l'application
\beq
\Ga^+(\si)\times \ba(Q) & \to &\Ga^+ (-\si)\times \ba(-Q)\\
			(\psi, A)& \mapsto& (\tau_S(\psi),\  A^*) 
\eeq
induit un isomophisme entre les solutions de 
$SW_h(\si,c)$ et celles de $SW_{-h}(-\si,c)$. On en déduit le résultat.
\epr
\section{ Calcul de l'invariant pour une variété symplectique}

Comme dans le cas classique, les travaux de Taubes
[Ta1] permettent de calculer
les invariants réels 
de Seiberg-Witten
pour certaines structures de $\spic$ sur une variété
symplectique réelle.

\subsection{Structures de $\spi$ sur une variété symplectique}
\subsubsection{Stucture canonique}
Soit $(X,\om,J)$ une variété symplectique munie d'une
structure presque complexe $J$ compatible avec $\om$, c'est à dire
 $g= \om(.,J.)$ est une métrique riemanienne. 
Gr\^ace à $J$, le groupe de transformation
$\SO$ de $P$ se réduit à $\U$ et permet de définir
 une structure naturelle de $\spic$ $\si_0$. 
En effet, le plongement $j: \U \to \spic(4)$ défini par
$$
j(A) = \left( \begin{array}{cc}
				\begin{array}{cc}
				1 & 0 \\	  
				0 & \det A  
				\end{array}		& 0\\
				0 & A \\
	\end{array} \right)
$$
permet de relever les fonctions de transitions $\gab \in \U$
à des fonctions de transitions $\tgab \in \spic(4)$ définissant
une structure de $\spic$ $\si_0$, de fibré déterminant
le fibré anticanonique $K^{-1}= \Lam^{0,2}$. Toutes les autres
structures de $\spic$ se déduisent de la canonique par la 
\begin{prop}Soit $\bl$ un fibré $\un$-principal sur $X$, et $L$ son
fibré en droites complexes associé. 
Alors  $\si_0 \otimes \bl$ est une structure de $\spic$ et on a les relations :
\beq 
S^+(L) &=& L \oplus K^{-1}\otimes L			\\
S^-(L) & = &\Lam^{0,1}\otimes L,\\
\det S^\pm(L) &= &  K^{-1}\otimes L^2.
\eeq
De plus, toute structure de $\spic$ est de la forme précédente.
\end{prop}

\subsubsection{Structures réelles}
Nous pouvons, dans ce cadre, expliciter quelques structures naturelles. 
Donnons d'abord  l'action de l'algèbre de Clifford $\Om(X,\Cc)$
sur $S^\pm(L)$. Si $e$ est un champ de vecteur
de $TX$ qu'on identifie avec  son
dual par la métrique, et $s\in \Om^{0,*}$, 
on a 
$$ e.s = \sqrt 2 (e^{0,1}\land s - e^{0,1}\angle \psi),$$
où $\angle$ est l'opérateur de contraction.

Ensuite soit $\nabla_{C}$ la connexion de Chern
sur $TX$. Cette connexion vérifie 
$\nabla_{C}g = 0$ ainsi que $\nabla_C J = 0$. Sa torsion 
est exactement le tenseur de Nijenhuis, si bien que $\nabla_C$
est la connexion de Levi-Cività dans le cas $\kah$.
$\nabla_C$ induit donc une connexion $A_0$ sur le fibré
anticanonique.

Supposons maintenant que la variété symplectique est réelle, 
c'est à dire qu'il existe une involution $c$ antisymplectique et 
$J$-antiholomorphe. 
Gr\^ace à la proposition précédente, on a  le 
\begin{lemme} 
Soit $\si = \si_0 \otimes \mathcal L$ une structure de $\spic$. 
Alors $\si$ est $c$-réelle si et seulement si
$\bl$ est $c$-réelle. En particulier, $\si_0$ est $c$-réelle, et
la connexion $\naz$ est $c$-r\'eelle. Enfin, l'opérateur
de Dirac $D_{A_0}$ vérifie :
\beq
D_{A_0} : \Om^0 \oplus \Om^{0,2} &\to& \Om^{0,1}\\
		(\al, \beta) & \mapsto &\bar \partial \alpha + \bar \partial^* \beta
\eeq
\end{lemme}
	\bpr
Si $\ka_\bl$ est une involution $c$-réelle sur $\bl$, 
l'application
\beq
\ka_S: \Ga^\pm(L) & \to &\Ga^\pm(L) \\
		\om \otimes s 	 & \mapsto &  \ka_\Om(\om)\otimes \ka_\bl(s)
\eeq
est une involution antilinéaire. 
 En remarquant que l'action de l'algèbre explicit\'ee plus haut commute
avec $\ka_\Om$, on trouve que
 $\ka_\Om (\om\cdot \psi) = \ka_\Om(\om)\cdot \ka(\psi).$

La connexion de Chern $\nabla_C$ est proportionnelle
à $\nabla + \frac{1}{2}J\nabla J$. On en déduit que tout comme
la forme de Levi-Cività, $\nabla_C$ commute avec $\ka_\Om$. 
En utilisant le lemme 3, on trouve que
 $\ka_{\mathcal A} (A_0) = \ka_S A_0 \ka_S = 
\ka_\Om A_0 \ka_\Om = \ka_\Om^2 A_0= A_0$, et donc que $A_0$ est réelle.
Enfin, on trouvera dans [Ni] tous les détails concernant 
la forme de Chern ainsi que la démonstration de l'égalité
entre $D_{A_0}$ et $\bar \partial  + \bar \partial^*$.
	\epr

		\subsection{Calcul de l'indice}
Rappelons que la dimension de l'espace
des solutions réelles est
est égale à $ -\dim a\bh^1 + \dim a\bhd + \frac{1}{2}\ind \di$.
Dans le cas d'une structure réelle sur une variété symplectique,
on a un calcul simple de la première partie de l'indice :
		\begin{lemme}
		Soit $(X,\om, J, c)$ une variété symplectique réelle compacte
		de dimension 4.
		On a 
		$ \dim a\bh^1 = \frac{1}{2}b_1$ et 
	 	$\dim a\bhd = \frac{1}{2}(1+ b_+)$.
		
		\end{lemme}
	\bpr
 Pour calculer
la dimension de $a\bh^1$, remarquons que l'application
$$\al \mapsto J^*\al = \star (\om \land \al)$$  
établit un isomorphisme entre  $a\Om^1$ et $r\Om^1$.
En effet, si $c^*\al = \al$, on a 
$c^*(J^*\al) = - J^*c^*\al = 
-J^*\al.$
Maintenant si $\al$ est harmonique, 
$J^*\al$ l'est aussi, car $\om$ est 
harmonique. On a donc la décomposition
$$ \bh ^1 = a\bh ^1 \oplus r \bh ^1.$$
Cela implique d'une part que $b_1$ 
est pair, d'autre part que 
$\dim a\bh^1 = \frac{1}{2} b_1.$

 Il nous reste
donc à calculer la dimension de $a\bhd$. 
Pour cela, rappelons la décomposition
des 2-formes autoduales :
$$ 
\Om^2_+ (X, \Rr) = 
\Om^0(X, \Rr)\ .\   \om \oplus 
\Om^{0,2},
$$
où $\om$ est la forme
symplectique et $  \Om^{0,2}$  l'espace (ici considéré comme
espace vectoriel sur $\Rr$) des (0,2)-formes complexes
sur $X$. 
Puisque $c$ est $J$-antiholomorphe,
l'application $\al \mapsto c^*\al$ envoie
$\Om^{0,2}$ sur $\Om^{2,0} = \Om^{0,2}$ en tant que espaces
réels. Si $\al \in \Om^{2,0}$ est  anti-$c$-réelle
et fermée, $J^*\al$ est $c$-réelle, fermée
et reste dans $\Om^{2,0}$. On a donc un isomorphisme
d'espace réels :
$$ 
   \Om^{2,0} \cap \ker d  = 
 a \Om^{2,0} \cap \ker d \oplus 
 r \Om^{2,0} \cap \ker d.
$$
L'égalité $c^*\om = - \om$ implique alors que $a\bhd \cap \Rr \om = \Rr \om$.
On a donc $\dim a\bhd = 1 + \frac{1}{2} \dim (\Om^{2,0}\cap \ker d).$
Sachant que $\dim \bhd = 1 + \dim (\Om^{2,0}\cap \ker d),$
on obtient $\dim a\bhd =  \frac{1}{2} (1 + b_+ ).$
	\epr
		\subsection{Calcul effectifs des invariants}
Nous pouvons énoncer maintenant le théorème principal 
concernant les variétés symplectiques :
\begin{theo} Soit $(X,\om,J,c)$ une variété symplectique compacte réelle,
vérifiant $b_+ >1$, et $K$ 
son fibré canonique. Alors $SW(\si_0,c) = \pm 1$ et 
$SW(\si_0\otimes K,c) = \pm 1$.
\end{theo}
\bpr
Remarquons d'abord que
la proposition 8 nous permet de ne démontrer que la première 
égalité.
La preuve est précisément celle du cas classique. Il 
faut simplement remarquer que les perturbations
choisies dans la preuve et la ou les solutions
\'el\'ementaires trouv\'ees sont $c$-réelles. 
 Nous rappelons uniquement les idées de la démonstrations qui 
suffisent \`a \'etablir le r\'esultat.

La perturbation $h$ choisie est $ih = F^+_{A_0} - \frac{i}{4}\rho^2 \om$, où $\rho$
est une constante strictement positive destinée à \^etre
très grande. La connexion
$A_0$ est $c$-r\'eelle, donc 
$\ka_\Om (F^+_{A_0}) = F^+_{A_0}$. Puisque de plus 
$c$ d\'efinit une structure r\'eelle
sur la vari\'et\'e symplectique, on a $\ka_\Om (\om) = - \om$, et
donc la perturbation choisie est $c$-réelle. 
Les équations étudiées sont donc, dans le cas plus général où $L$ n'est
pas forcément triviale:
\beq
(\psi, A) & \in & r\Ga^+\times r\ba\\
\di \psi& = &0\\
F^+_A &= &q (\psi)+ F^+_{A_0}-\frac{i}{4}\rho^2 \ \om,
\eeq
Sur le fibré déterminant $Q= K^{-1}\otimes L^2$, une connexion $A$ peut s'écrire
$A= A_0\otimes B^2$, où $B$ est une connexion sur $L$.  On 
utilise $B$ comme nouvelle variable au lieu de $A$. 
Les équations dans le cas symplectique presque complexe se transforment en :
\beq
((\al, \beta),B ) &\in  &r\Om^0 (L) \oplus r \Om^{0,2}(L) \times r\ba(L)\\
\del_B \al + \del_B^*\beta &=& 0 \\
(F_B^+)^{(1,1)} = \frac{i}{8} (|\al|^2 - |\beta|^2 - i\rho^2) \om \ & \text{et}& \ 
F_B^{(0,2)} = \frac{\bar \al \beta}{4}.
 \eeq
Si comme dans l'énoncé 
$L$ est triviale, on a une solution
manifestement invariante par $\ka_S$ et $\ka_\ba$ de ces \'equations :
$$((\al, \beta), B) = ((\rho, 0), d).$$
On d\'emontre (cf. [Ta1]) que pour $\rho$ assez grand, 
c'est l'unique solution
modulo l'action de $\mathcal G$.
Or deux solutions $c$-r\'eelles $\mathcal G$-\'equivalentes
sont n\'ecessairement $r\mathcal G$-\'equivalentes,
ce qui finit de d\'emontrer le th\'eor\`eme.
\epr

Nous parasitons le difficile théorème de Taubes pour démontrer :
\begin{theo}
Soit $(X,\om,J,c)$ une variété symplectique compacte réelle, 
vérifiant $b_+ >1$ et $L$ 
un fibré en droite $c$-réel sur $X$.
Si $SW(\si_0\otimes \bl,c) \not= 0$, 
et $\bl \not=0$, alors il existe une
courbe $J$-holomorphe réelle (a priori singuli\`ere et non connexe) représentant
le fibré $\bl$. 
\end{theo}
\bpr
Les arguments de Taubes développés dans [Ta1] s'appliquent
entièrement à notre situation. Il faut simplement
traduire géométriquement le fait qu'on a des solutions
réelles.

Si l'invariant de Seiberg-Witten r\'eel $SW(\si,c)$ est non nul, 
alors pour toute
perturbation d\'efinie comme ci-dessus, c'est \`a dire pour
toute suite $\rho_n $ tendant vers 
l'infini, il existe une solution r\'elle
aux \'equations $((\al_n, \beta_n),B_n)$.
Taubes montre
d'une part que la courbure
$F_{B_n}$ converge au sens des courants 
vers une courbe $J$-holomorphe $C$ a priori singulière et
non connexe. D'autre part 
il montre que le lieu  des z\'eros de $\al_n$ tend au sens
de Hausdorff vers $C$. Les $\al_n$ étant des sections
du fibré $\bl$, $C$ est Poincaré duale à $\bl$.

On sait de plus que le spineur $(\al_n, \beta_n)$
est invariant par $\ka_S= \ka_\Om\otimes \ka_\bl$. En particulier  
$\ka_\bl \al_n = \al_n$, ce qui implique que pour tout $n$,
le lieu des z\'eros de $\al_n$ 
est invariant par $c$.  La limite au sens de Hausdorff, en l'occurence
$C$, l'est aussi, et 
le th\'eor\`eme est d\'emontr\'e.
\epr

\begin{coro}
Soit $(X,\om,J,c)$ une variété symplectique compacte réelle, 
vérifiant $b_+ >1$. Alors 
le fibré canonique est représenté par une courbe réelle.
\end{coro}

\textsc{D. Gayet: Mathématique, B\^atiment 425, Universit\'e de Paris Sud,
 91405 Orsay Cedex France.} 

E-mail : damien.gayet@math.u-psud.fr
\end{document}